\documentclass[12pt,a4paper]{article}

\usepackage{amsmath,amsthm,amsfonts}

\usepackage{setspace}

\usepackage{graphicx}

\def\back{\noindent\kern-5mm} 

\def\de{{\rm d}}

\def \Label{\label}
\def \bE{{\mathbf E}}
\def \bP{{\mathbf P}}
\def \bR{{\mathbb R}}

\def \de{{\rm d}}

\def \mS{{\mathcal S}}
\def \mP{{\mathcal P}}
\def \mR{{\mathcal R}}
\def \mE{{\mathcal E}}
\def \mM{{\mathcal M}}

\def \zs#1{_{\lower 3pt \hbox{$\scriptstyle#1$}}}

\def\beq{\begin{equation}}
\def\endeq{\end{equation}}
\def\beqn{\begin{eqnarray*}}
\def\endeqn{\end{eqnarray*}}

\newtheorem{theorem}{Theorem}

\def\thetitle{Goodness of fit test for ergodic diffusion processes \footnote{This work has been partially supported by the local grant sponsored by University of Bergamo: {\it Theoretical and computational problems in statistics for
continuously and discretely observed diffusion processes} and by MIUR 2004 Grant. } }

\title\thetitle
\author{
Ilia Negri\footnote{%
Corresponding author: University of Bergamo, Viale Marconi, 5, 24044, Dalmine (BG), Italy.
{\tt ilia.negri@unibg.it}.} 
\and Yoichi Nishiyama\footnote{%
The Institute of Statistical Mathematics, 4-6-7, Minami-Azabu, Minato-ku, Tokyo, 106-8569, Japan.
{\tt nisiyama@ism.ac.jp}}
}

\begin{document}


\maketitle

\begin{abstract}
A goodness of fit test for the drift coefficient of an ergodic diffusion process is presented. The test is based on the  score marked empirical process.  The weak convergence of the  proposed test statistic is studied under the null hypotheses and it is  proved that the limit process is a continuous Gaussian process. The structure of  its covariance function allows to calculate the limit distribution and it turns out that it is a function of a standard Brownian motion and so exact reject regions can be constructed. The proposed test is asymptotically distribution free
and it is consistent under any simple fixed alternative.
\end{abstract}

\noindent
{\bf Key words:} consistent test, empirical process, asymptotically distribution free tests.
\par

\noindent
{\bf 2000 MSC:} 60G10; 60G35; 62M02.
\par

\vfill
\eject

\doublespacing

\section{Introduction}
Let $X$ be an ergodic diffusion process on $\bR$, solution of a stochastic differential equation, that is a strong Markov process with continuous sample paths  which satisfies
$\de X_t=  S(X_t) \de t+\sigma(X_t)\de W_t$, for $t>0$, with some random initial value $X_0$. We consider the nonparametric problem of testing the unknown drift coefficient $S$ on the basis of the continuous observation $X^T=\{X_t: 0\leq t\leq T\}$. We present a goodness of fit test for the drift coefficient of a diffusion process based on the statistic $\sup_x |V_T(x)|$ where 
$$
V_T(x)=\frac{1}{\sqrt{T}} \int_0^T{ \bf 1}\zs{(-\infty,x]}(X_t)(\de X_t- S_0(X_t)\de t).
$$
As usual  ${\bf 1}_A$ denote the indicator function on a set $A$. Following Koul and Stute (1999), we call $V_T$ the {\it score marked empirical process}. We prove that the test based on the statistic $\sup_x |V_T(x)|$ is asymptotically distribution free under each simple null hypothesis  $S=S_0$
and it is consistent under any simple fixed alternative $S=S_1\neq S_0$.  

Despite the fact of  their importance in applications, few works are devoted to the  goodness of fit test for diffusions up to our knowledge. So the construction of goodness of fit tests for such kind of model is very important and needs very detailed studies.  
Kutoyants (2004) discusses some possibilities of the construction of such tests. In particular,  he considers the Kolmogorov-Smirnov statistics
$\Delta_T(X^T)= \sup_x \sqrt{T}|\hat F_T(x)-F_{S_0}(x)|$. The goodness of fit test based on this statistics is asymptotically consistent and the asymptotic distribution under the null hypothesis follows from the weak convergence of the empirical process to a suitable Gaussian process (Negri, 1998 and Van der Vaart and Van Zanten, 2005). However, due to the structure of the covariance of the limit process, the Kolmogorov-Smirnov statistics  is not asymptotically distribution free in diffusion process models. The statistic proposed here does not have this problem.

We study the weak convergence of the  score marked empirical process under the null hypotheses and we prove that the limit process is a continuous Gaussian process. The structure of  its covariance function allows us to calculate the limit distribution of the proposed statistic. It turns out that it is a functional of a standard Brownian motion with known distribution and so we can construct exact reject regions.  

Koul and Stute (1999) proposed such kind of statistics based on a class of empirical process constructed on certain  residuals to check some parametric models for time series. They  studied their  large sample behavior under the null hypotheses and present a martingale transformation of the underlying process that makes tests based on it asymptotically distribution free. Some considerations on consistency have also been done. For the same model studied here the problem of testing different parametric form of the drift coefficient is well developed (see  for example Lin'kov, 1981,  Kutoyants, 2004 and references therein). 

The work is organized as follows. In the next section we present the model, its properties and some general conditions and assumptions used through all the text. In Section  \ref{limit} we prove the weak convergence of the proposed statistic in a more general context. This result is interesting by itself. Section \ref{test} is devoted to the presentation of the test problem and to the study of the proposed statistic under the null hypothesis.  Finally in Section \ref{asy} we study the behavior of the statistics under the alternative hypothesis and we prove that the test is consistent.

\section{Preliminaries}
\Label{sec2}
Given a general stochastic basis, that is, a probability space  $(\Omega, {\cal A}, \bP)$ and a filtration 
$\{{\cal A}_t\}_{t\geq0}$ of ${\cal A}$, let us consider a one dimensional diffusion process solution of the following stochastic differential equation
\begin{equation}
\begin{cases}
\de X_t = S(X_t) \de t+\sigma(X_t)\de W_t\\ 
X_0=\xi,  
\end{cases}
\Label{eds}
\end{equation}
where $\{W_t:\ t\geq 0\}$ is a standard Wiener process, and the initial value $X_0=\xi$ is independent of $W_t$, $t\geq 0$.
The drift coefficient  $S$ will be supposed unknown to the observer
and the diffusion coefficient $\sigma^2$
will be a known positive function. 
Let us introduce the following condition.

$\mE\mS$. {\em The function $S$ is locally bounded, the function $\sigma^2$ is continuous and bounded and for some constant $A>0$, the condition
$ xS(x)+\sigma(x)^2\leq A(1+x^2)$, $x\in \bR$, holds.}

Under condition $\mE\mS$ the equation \eqref{eds} has an unique weak solution (see Durrett, 1996, p. 210). 
The {\em scale function} of a diffusion process solution of the stochastic differential equation \eqref{eds}
is defined by
$$
p(x)=\int_{0}^{x}\exp\left\{-2\int _0^y{S(v)\over \sigma ^2(v) }\; {\rm d}v
\right\}  {\rm d}y.
$$ 
The {\em speed measure} of the diffusion process \eqref{eds} is defined by $m_S(\de x)= \frac{1}{\sigma(x)^2 p'(x)}\de x$. 
Let us introduce the following condition:

$\mR\mP$. {\em The scale function is such that
$$
\lim_{x\to \pm \infty} p_(x)= \pm \infty
$$
and the speed measure $m_S$ is finite. 
}

If the condition $\mR\mP$ is satisfied then the process  $\{X_t: t\geq 0\}$, weak solution of \eqref{eds}, has the ergodic property (see for example Gikhman and Skorohod, 1972 or Durrett, 1996), that is, there exists an unique invariant probability measure $\mu_S$ such that for every
measurable function $g\in {\cal L}_1(\mu_S)$ we have with  probability one,
$$
 \lim_{T\to\infty}\frac1T \int_0^T g(X_t)\de t =
\int_\bR g(z) \mu_S( \de z).
$$
Moreover the invariant measure $\mu_S$ has a density given by 
$$
f_S(y)=\frac{1}{m_S(\bR)\sigma(y)^2}
\exp\left\{2\int_0^y \frac{S(v)}{\sigma(v)^2}\de v \right\},
$$
where
 $$
m_S(\bR)= \int_{-\infty}^{\infty}\frac{1}{\sigma(x)^{2}}  \exp\left\{2\int
_0^x{S(v)\over \sigma ^2(v) }\; {\rm d}v \right\} \; {\rm d}x 
$$
is finite.
\section{A limit theorem}
\label{limit}
In this section we present a theorem on the weak convergence of a stochastic process that is interesting by itself. Let us consider a diffusion process $X=\{X_t :t\geq 0\}$ on an open interval $I\subseteq \bR$, that is a strong Markow process with continuous sample paths taking values on $I$, not necessarily solution of a differential stochastic equation of type \eqref{eds}. Assume that  $X$ is regular, which implies that the scale function $p$ and the speed measure $m$ of the diffusion are well defined (see Rogers and Williams, 2000). Under the assumption that the speed measure  $m$ is finite, and denoting by $\mu$ the normalized speed measure, $\mu=\frac{m}{m(I)}$, it follows that the diffusion process $X$ is positive recurrent and it has the ergodic property, with $\mu$ as invariant measure.  
For every $x\in I$ and $t\geq0$, the local time for the diffusion $X$ in the point $x$ at time $t$ is denoted $L_t^X(x)$. The process $\{L^X_t(x) : x\in I, t\geq 0\}$, is countinuous and increasing in $t$ and  continuous and nonnegative in $x$.   
The main theorem for local time is the {\it occupation time} formula. For a diffusion process it can be written as 
\begin{equation}
\int_0^t h(X_s) ds =\int_I L_t^X(x)h(x)p'(x)m(dx).
\label{of}
\end{equation}
for every  measurable function $h:I \to \bR$. Let us denote the diffusion local time for $X$ with respect to the speed measure (see Van der Vaart and Van Zanten, 2005 and It\^o and McKean, 1965) as $l_t^X(x)=L_t^X(x)p'(x)$. The random function $x\to l_t^X(x)$ can be chosen continuos and has compact support.
If  the measure $m$ is finite  then  it holds that  
\begin{equation}
\frac{1}{t}\sup_{x\in I}l_t^X(x)=O_\bP(1).
\label{VVZ}
\end{equation}
See theorem 4.2 of Van der Vaart and Van Zanten (2005) and also Van Zanten (2003). 

Let a standard Wiener process $W$ be given on the same stochastic basis where $X$ is defined.  Let us consider the process $M=\{ M_t(\psi): t\geq 0, \psi \in {\cal F}\}$ defined by 
$$
M_t(\psi)=\int_0^t \psi(X_s)\de W_s,
$$
where $\psi$ belongs to a countable  class ${\cal F}$ of elements of  ${\cal L}^2(I, m(\de x))$.  To measure the distance between functions in  ${\cal F}$ we use the semimetric $\rho$ 
$$
\rho(\psi, \varphi)= \sqrt{\int_I |\psi(y)- \varphi(y)|^2 m(\de y)}.
$$
For every $\psi\in {\cal F}$ the process $M(\psi)=\{M_t(\psi): t\geq 0\}$ is a continuous local martingale. Following Nishiyama (1999), definition 2.1, the quadratic  $\rho$-modulus for the process $V=\{V_t(\psi) : t\geq 0, \psi \in {\cal F}\}$, where $V_t(\psi)= \frac{1}{\sqrt{t}}M_t(\psi)$,  is defined as
$$
||V||_{\rho, t}=\sup_{\rho(\psi, \varphi)>0}\frac{\sqrt{\frac{1}{t}\left\langle M(\psi)-M(\varphi)\right \rangle_t}}{\rho(\psi, \varphi)}.
$$ 
Here $\langle M\rangle=\{ \langle M\rangle_t:t\geq 0\}$ denotes the quadratic variation process
of  a continuous local  martingales $M$. 

Let us denote with $N (\epsilon, {\cal F}, \rho)$ the smallest number of closed balls, with $\rho$-radius $\epsilon>0$, which cover the set ${\cal F}$. 
\begin{theorem}
Let ${\cal F} \subset {\cal L}^2(I, m(\de x))$ be countable. Then for all $\delta$, $K$ and  $t>0$ it holds that
$$
\bE \sup_{\rho(\psi, \varphi)<\delta} | V_t(\psi)- V_t(\varphi)| {\bf 1}\zs{\{ \xi_t \leq K\}} \leq cK \int_0^\delta \sqrt{ \log N (\epsilon, {\cal F}, \rho) } \de \epsilon,
$$
where $c>0$ is an universal constant and $\{\xi_t: t>0\} $ is a stochastic process which satisfies $\xi_t =O_\bP(1)$, as $t$ goes to infinity. 
\label{teo1}
\end{theorem}
\begin{proof}
We have
$$
\left< V(\psi)\right>_t= \frac{1}{t}\int_0^t  \psi(X_s)^2  \de s
$$
and from the occupation formula \eqref{of} for the diffusion local time of a diffusion process we can write
$$
\frac{1}{t}\int_0^t  \psi(X_s)^2\de s =\frac{1}{t} \int_I \psi(x)^2 l_T^X(x)m(\de x).
$$
So we have
$$
\| V\|_{\rho, t}= \sup_{\rho(\psi, \varphi)>0} \frac{\sqrt{\frac{1}{t}\int_I (\psi(x)- \varphi(x))^2 l_T^X(x)m(\de x)}}{\rho(\psi, \varphi)}\leq \sqrt{\frac{1}{t} \sup_{x\in I}  l_T^X(x)}.
$$
Now recalling \eqref{VVZ} 
the result follows from the theorem 2.3 in Nishiyama (1999) if we pose $\xi_t= \sqrt{\frac{1}{t}\sup_{x\in I}  l_t^X(x)}$.  
\end{proof} 
Let us denote with $\ell^{\infty}({\cal F})$ the space of bounded functions ${\cal Z}: {\cal F}\to \bR$ equipped with the uniform norm $||{\cal Z}||_{\infty}=\sup_{\varphi \in {\cal F}} |{\cal Z}(\varphi)|$.

On the space $ \ell^{\infty}({\cal F})$ we introduce the Gaussian process $\{\Gamma(\psi): \psi \in {\cal F}\}$  with mean zero and covariance function given by
$$
g(\psi,\varphi)= \int_I \psi(z)\varphi(z)\mu(\de z).
$$

\begin{theorem}
Let  ${\cal F} \subset {\cal L}^2(I, m(\de x))$ be countable and $ \int_0^1 \sqrt{ \log N (\epsilon, {\cal F}, \rho) } \de \epsilon$  be finite. Then the family of stochastic maps $\{V_t(\psi): \psi \in {\cal F}\}$ weakly converges on the space $\ell^{\infty}({\cal F})$, as $t$ goes to infinity, to the Gaussian process  $\{\Gamma(\psi): \psi \in {\cal F}\}$.
\label{wcV}
\end{theorem}
\begin{proof}
The convergence of the finite dimensional laws $(V_t(\psi_1), \ldots V_t(\psi_k))$ to the law of $(\Gamma(\psi_1), \ldots,\Gamma(\psi_k))$, for every finite $k$ follows from the central limit theorem for stochastic integrals (see Kutoyants, 2004). The tightness follows from  Theorem \ref{teo1}, assuming  that$ \int_0^1 \sqrt{ \log N (\epsilon, {\cal F}, \rho) } \de \epsilon$  is finite (see also Nishiyama, 2000).
\end{proof}

\section{Goodness of fit test}
\label{test}
Let us introduce our testing problem. Suppose that we observe the process $\{X_t : 0\leq t\leq T\}$, solution of the stochastic differential equation \eqref{eds} and we wish to test the two different simple hypotheses 
$$
\begin{array}
[c]{l}%
H_{0}: S=S_0\\
H_{1}: S =S_1,%
\end{array}
$$
where $S_1\neq S_0$ means 
$$
\sup_{x\in \bR} |S_1(x)-S_0(x)|>0.
$$ 
We suppose that $S_0$ and $S_1$ belong to the class ${\cal S}_\sigma$ defined
for a fixed function $\sigma$ as
$$
{\cal S }_\sigma=\left\{  S : \text{conditions } \mE\mS \text{ and } \mR\mP \text{ are fulfilled} \right\}.
$$ 
For every $x \in \bR$, let us introduce the {\em score marked empirical process}
\begin{eqnarray*}
V_T(x)& = &\frac{1}{\sqrt{T}} \int_0^T {\bf 1}\zs{(-\infty,x]}(X_t)(\de X_t- S_0(X_t)\de t)\\
&=& \frac{1}{\sqrt{T}} \int_0^T {\bf 1}\zs{(-\infty,x]}(X_t)\sigma(X_t)\de W_t.
\end{eqnarray*}
The process $\{V_T(x): x\in \bR\}$ takes values in $C_B(\bR)$, the space of the continuous bounded function on $\bR$. Let us introduce in this space the $\sigma$-algebra of borel set ${\cal B}$ generated by the open sets of $C_B(\bR)$ induced by the norm $||f||=\sup_{x\in \bR}|f(x)|$ for every $f\in C_B(\bR)$. 

Let us introduce our test procedure. Fix a number $\varepsilon \in (0,1)$ and let us consider the class of {\em asymptotic test of level} $1-\varepsilon$ or {\em size} $\varepsilon$. Given any statistical decision function $\phi_T=\phi_T(X^T)$, 
the expected value of  $\phi_T(X^T)$ is the probability to reject $H_0$ having the observation $X^T=\{X_t : 0\leq t\leq T\}$.  Let us denote by $\bE_S^T$  the
mathematical expectation with respect to 
the measures $\bP_S^T$ induced by the process $\{X_t :  0\leq t\leq T\}$ in the space $C[0,T]$ (the space of all the continuos functions on $[0,T]$). We define the class of all the test of asymptotic level $1-\varepsilon$ as
$$
{\cal K}_\varepsilon=\left\{ \phi_T : \limsup_{T\to +\infty} \bE_{S_0}^T\phi_T(X^T)\leq\varepsilon \right\}.
$$
The power function of the test based on $\phi_T$ is the probability of the true decision under $H_1$, and is given by
$$
\beta_t(\phi_T)=\bE_{S_1}^T\phi_T(X^T).
$$ 
A test procedure is consistent if 
$$
\lim_{T\to +\infty} \bE_{S_1}^T\phi_T(X^T) =1.
$$
Let us introduce in the space   $(C_B(\bR), {\cal B})$ the Gaussian process $\{ \Gamma(x): x \in \bR\}$ 
with mean zero and covariance function given by $G(x,y)=g^2_{S_0}(x \wedge y)$, where
$$
g^2_{S_0}(z)=\int_{-\infty}^{+\infty} {\bf 1}\zs{(-\infty,z]}(y)\sigma(y)^2 f_{S_0}(y)\de y 
$$
Note that $\sup_z g^2_{S_0}(z)$ is finite if $\sigma^2$ is bounded. 
Since the function $g_{S_0}$ is nondecreasing and nonnegative the limit process $\{ \Gamma(x): x \in \bR\}$ admits the following representation in distribution
\begin{equation}
\Gamma(x)= B(g_{S_0}(x))
\label{brownian}
\end{equation}
for every $x\in \bR$, where $B$ denote a standard Brownian motion on the positive real line.  
The weak convergence of the process $\{V_T(x):x\in \bR\}$ to the process $\{ \Gamma(x): x \in \bR\}$  in the space 
$(C_B(\bR), {\cal B})$ is immediate if we apply Theorem \ref{wcV} to the functions $\psi(y)={\bf 1}_{(-\infty, x]}(y)\sigma(y)$.  Here we remark that the class ${\cal F}$ in Theorem \ref{wcV} has to be countable. However  the process $\{V_T(x):x\in \bR\}$ is continuous in $x$, so in the current situation we can consider such class of functions.
This result on the weak convergence of process $V_T$,  the continuous mapping theorem and the representation \eqref{brownian} yield the following relations in distribution
$$
\sup_{x\in \bR} |V_T(x)| \Rightarrow \sup_{x \in \bR} |\Gamma(x)|= \sup_{0\leq t\leq g^2_{S_0}(\infty)} |B(t)|=  g_{S_0}(\infty)\sup_{0\leq t \leq 1}|B(t)|
$$
where $g^2_{S_0}(\infty)= \int_\bR \sigma(y)^2  f_{S_0}(y) \de y$ is finite, and $\Rightarrow$ denote the weak convergence as $T$ goes to infinity.
We will consider the following statistical decision function
$$
\phi_T^*= {\bf 1} \zs{\left\{\frac{1}{g_{S_0}(\infty)}\sup_{x \in \bR}| V_T(x)|>c_\varepsilon \right\}}
$$
where  the {\it critical value} $c_\varepsilon$ 
is defined by 
$$
\bP\left( \sup_{0\leq t\leq 1} |B(t)|>c_\varepsilon\right)=\varepsilon. 
$$
So we have proved that $\phi_T^*\in {\cal K}$ and that the test is asymptotically distribution free. 
To be the introduced statistical procedure useful we have to study the asymptotic properties  of the statistics $\sup_x |V_T(x)|$ under  the alternative hypotheses. This is done in the next section.

\section{Consistency of the test}
\label{asy}
This section is devoted to the study  of the asymptotic beaviour of the test statitistic  $\sup_{x \in \bR} |V_T(x)|$ under  the alternative hypotheses. We have shown in the previous section  that the proposed statistics is asymptotically distribution free. Now we prove that the proposed test procedure is also consistent under very general condition on the model.    
Let us introduce the following condition.

${\cal C}$: For some $x\in \bR$ it holds
$$
\int\limits_{-\infty}^{+\infty} {\bf 1}\zs{(-\infty,x]}(y)(S_0(y)-S_1(y))f_{S_1}(y)\de y\ne 0.
$$ 
\begin{theorem}
\label{MT}
Let $S_0$ and $S_1$ belong to ${\cal S}_\sigma$ and the condition ${\cal C}$ be satisfied. Then the test based on the statistical decision function
$$
\phi_T^*= {\bf 1} \zs{\left\{\frac{1}{g_{S_0}(\infty)}\sup_{x \in \bR}| V_T(x)|>c_\varepsilon \right\}}
$$
is consistent. 
\end{theorem}
\begin{proof}
To prove the consistency it is enough to show that, under $H_1$
$$
\bP\left( \lim_{T\to +\infty}\sup_{x\in \bR}|V_T(x)|= +\infty \right)=1.
$$
We can write
$$
\sup_{x\in \bR}|V_T(x)|\geq \sqrt{T}\sup_{x\in \bR}|A_T(x)| - \sup_{x\in \bR}| V^1_T(x)|,
$$
where $V^1_T(x)$ and $A_T(x)$ are given as follows.
Under $H_1$, by Theorem \ref{wcV} the process 
$$
V^1_T(x)=\frac{1}{\sqrt{T}} \int_0^T {\bf 1}\zs{(-\infty,x]}(X_t)(\de X_t- S_1(X_t)\de t)
$$
weakly converges to the corresponding Gaussian process so the limit process is tight. 
On the other hand,  
$$
A_T(x) =\frac{1}{T}\int_0^T
{\bf 1}\zs{(-\infty,x]}(X_t)(S_0(X_t)-S_1(X_t))\de t
$$
converges a.s. uniformly in $x$ to 
$$
A(x)=\int\limits_{-\infty}^{+\infty} {\bf 1}\zs{(-\infty,x]}(y)(S_0(y)-S_1(y))f_{S_1}(y)\de y.
$$
If the condition ${\cal A}$ is satisfied we have 
$$
\lim_{T\to +\infty}\sqrt{T}\sup_{x\in \bR}|A_T(x)|= +\infty \quad {\rm a. s.}
$$
and the test is consistent.  
\end{proof}




\section*{References}
{\leftskip5mm

\back  Durrett, R.  (1996). {\sl Stochastic Calculus: A Practical Introduction}, CRC Press, Boca Raton.

\back  Gikhman, I.I., Skorohod, A.V. (1972). {\sl Stochastic
        Differential Equations}, Springer--Verlag, New York.

\back Koul, H.L., Stute, W. (1999). Nonparametric model checks for times series, {\sl The Annals of Statistics}, {\bf 27}, 1, 204-236.

\back Kutoyants, Y.A. (2004). {\sl Statistical Inference for Ergodic Diffusion Processes}, Springer, New York.

\back It\^o, K., McKean, H.P. (1965). {\it Diffusion Processes and Their Sample Paths},
Springer-Verlag, Berlin.

\back Lin'kov, Y.N. (1981). On the asymptotic power of a statistical test for diffusion type processes. (Russian) {\it Theory of random processes}, {\bf 9},  61--71. 
        
\back Negri, I. (1998). Stationary distribution
function estimation for ergodic diffusion process, 
        {\sl Statistical Inference for Stochastic Processes}, {\bf 1}, 61--84.

\back Nishiyama, Y. (1999). A maximal inequality for continuous martingales and M-estimators in a Gaussian white noise model, {\it The Annals of Statistics}, {\bf 27}, 2, 675-696.  

\back Nishiyama, Y. (2000). {\em Entropy Methods for Martingales}, CWI Tract {\bf 128}, Centrum voor Wiskunde en Informatica, Amsterdam.

\back Rogers, L.C.G., Williams, D. (2000). {\it Diffusion, Markov Processes and Martingales}, Vol. II, Cambridge University Press. 

\back Van der Vaart, A.W.,  Van Zanten, H.  (2005). Donsker theorems for diffusion: necessary and sufficient conditions, {\it The Annals of Probability}, {\bf 33}, No. 4, 1422-1451. 

\back Van Zanten H.  (2003). On empirical processes for ergodic diffusions and rates of convergence of M-estimators, {\it Scandinavian Journal of Statistcs}, {\bf 30},  443-458. 

}

\end{document}